\documentclass[oneside]{amsart}
\usepackage{amssymb,cite}
\numberwithin{equation}{section}

\usepackage[a4paper,
includeheadfoot]{geometry}

\usepackage{latexsym}
\usepackage{textcomp}
\usepackage{amsmath}
\usepackage{amsfonts}
\usepackage{amssymb}
\usepackage{mathrsfs}
\usepackage{enumerate}

\renewcommand{\a }{\alpha }
\renewcommand{\b }{\beta }
\renewcommand{\d}{\delta }

\newcommand{\pa}{{\partial}}
\newcommand{\D }{\Delta }

\newcommand{\e }{\varepsilon }
\newcommand{\g }{\gamma}

\newcommand{\G }{\Gamma}
\renewcommand{\l }{\lambda }

\newcommand{\n }{\nabla }
\newcommand{\vp }{\varphi }
\renewcommand{\phi}{\varphi}

\newcommand{\s }{\sigma }

\renewcommand{\t }{\tau }

\renewcommand{\th }{\theta }

\renewcommand{\O }{\Omega }

\newcommand{\ov}{\overline}

\newcommand{\be}{\begin{equation}}
\newcommand{\ee}{\end{equation}}

\newcommand{\R}{\mathbb{R}}
\newcommand{\N}{\mathbb{N}}
\renewcommand{\S}{\mathbb{S}}

\newcommand{\de}{\partial}

\newcommand{\ti}{\widetilde}

\renewcommand{\k}{\kappa}

\newcommand{\calD }{\mathcal{D}}

\def\dlim{\displaystyle\lim}

\newtheorem{Theorem}{Theorem}[section]

\newtheorem{Lemma}[Theorem]{Lemma}

\newtheorem{Corollary}[Theorem]{Corollary}
\newtheorem{Remark}[Theorem]{Remark}

\def\proof{\noindent{{\bf Proof. }}}
\def\square{\vbox{
    \hrule height .4pt
    \hbox{\vrule width .4pt height 7pt \kern 7pt
       \vrule width .4pt}
    \hrule height .4pt }}

\def\square{\vbox{
    \hrule height .4pt
    \hbox{\vrule width .4pt height 7pt \kern 7pt
       \vrule width .4pt}
    \hrule height .4pt }}

\def\QED{\hfill {$\square$}\goodbreak \medskip}

\def\R{{\mathbb R}}

\def\S{{\mathbb S}}

\def\div{{\rm div}}

\newcommand{\Ds}{(-\Delta)^s}
\newcommand{\Dsm}{(-\Delta+m^2)^s}
\newcommand{\RNp}{\R^{N+1}_+}

\newcommand{\SN}{{\mathbb S}^{N-1}}

\newcommand{\dive }{\mathop{\rm div}}

\begin{document}
\title[Relativistic Schr\"odinger operators with a singular
potential]{Sharp essential self-adjointness   of relativistic Schr\"odinger operators with a singular
  potential} \author[Mouhamed Moustapha Fall \and Veronica
Felli]{Mouhamed Moustapha Fall \and Veronica Felli }
\address{\hbox{\parbox{5.7in}{\medskip\noindent
      M.M. Fall\\
      African Institute for Mathematical Sciences (A.I.M.S.) of Senegal,\\
      KM 2, Route de Joal,\\
      B.P. 1418.
      Mbour, S\'en\'egal. \\[2pt]
      {\em{E-mail address: }}{\tt mouhamed.m.fall@aims-senegal.org.}\\[5pt]
      V. Felli\\
      Universit\`a di Milano
      Bicocca,\\
      Dipartimento di Ma\-t\-ema\-ti\-ca e Applicazioni, \\
      Via Cozzi
      55, 20125 Milano, Italy. \\[2pt]
      \em{E-mail address: }{\tt veronica.felli@unimib.it.}}}}

\thanks{M. M. Fall is supported by the Alexander von Humboldt
  foundation.  V. Felli was partially supported by the P.R.I.N. 2012 grant
  ``Variational and perturbative aspects of nonlinear differential
  problems''.
  \\
  \indent 2010 {\it Mathematics Subject Classification.} 35R11, 35B40,
  47B25, 35J75.\\
  \indent {\it Keywords.} Fractional elliptic equations, essential
  self-adjointness, singular homogeneous potentials.}

\date{May 8, 2014}

 \begin{abstract}
   \noindent
  This paper is devoted to  the study of essential self-adjointness of
  a relativistic Schr\"{o}dinger operator
  with a singular homogeneous potential.
  From an explicit condition on the coefficient of the singular term,
  we  provide a sufficient and   necessary  condition for    essential self-adjointness.
 \end{abstract}

\maketitle

\section{Introduction}
 The purpose of the present paper is to provide sharp  essential self-adjointness
of the Hamiltonian
\begin{equation}\label{eq:Ham}
H(p,x):=(p^2+m^2)^s-\frac{a(\frac{x}{|x|})}{|x|^{2s}},\quad
x\in\R^N,
\end{equation}
with $a:\SN\to\R$, $s\in(0,1)$, $m\geq0$, $N>2s$.

A symmetric densely-defined operator in a Hilbert space is said to be
\emph{essentially self-adjoint} if it has a unique self-adjoint
extension. We recall that if a symmetric operator $A: D(A)\to E$,
with $D(A)$ dense in the Hilbert space $E$,  is  strictly
positive, i.e. if $(Au,u)_E\geq c(u,u)_E$ for all $u\in D(A)$ and some
$c>0$, then $A$ is essentially self-adjoint if and only if its range
is dense  in $E$, see e.g. \cite[Theorem X.26]{reedsimon}.

In $3$-space dimension, the quantum
mechanics of a spin zero relativistic particle of charge $e$ and mass
$m$ in the Coulomb field of an infinitely heavy nucleus of charge $Z$
is described by the Hamiltonian $H(p,x)=(p^2+m^2)^{1/2}-Z e^2|x|^{-1}$, see e.g. \cite{Herbst,lieb}.
From \cite{Kato, Herbst} it is known that $(p^2+m^2)^{1/2}-Z
e^2|x|^{-1}$ is semi-positive definite if $Ze^2\leq 2/\pi$ and,
moreover, it is   essentially-self adjoint if $Ze^2\leq 1/2$. As a
particular case of the main result of the
present paper, we will see that if $Ze^2> 1/2$
then  $(p^2+m^2)^{1/2}-Z
e^2|x|^{-1}$ is not  essentially-self adjoint.

The essential self-adjointness of the operator $A
=H(i\nabla,x)=\Dsm-a(x/|x|) |x|^{-2s} $ implies uniqueness of the
quantum dynamics defined by $A$. Next, once we know that an operator
is not essentially self-adjoint, the choice of its extension to
generate the quantum dynamics is dictated by the physics problem, see
\cite{reedsimon} for more explanations.  Another application of
essential self-adjointness is in probability. Indeed, in general, $A$
could have several self-adjoint extensions $A'$, yielding Markov
processes with transition semigroups $p_t = e^{-tA'}$. The essential
self-adjointness of $A$ implies that there is only one self-adjoint
extension $A_F$: the \emph{Friedrichs extension}.  Hence, in case of
essential self-adjointness, we have a unique such semigroup and thus a
unique Markov process with generator $A_F$.

Let  ${\mathbb S}^N$ be the unit $N$-dimensional sphere  and
\[
{\mathbb S}^{N}_+=\{(\theta_1,\theta_2,\dots,\theta_N)\in
{\mathbb S}^N:\, \theta_1>0\}.
\]
We will denote by $dS$ (respectively  $dS'$) the volume
element on $N$-dimensional (respectively $(N-1)$-dimensional)
spheres and define
$H^{1}({\mathbb S}^{N}_+;\theta_1^{1-2s})$ as the completion of
$C^\infty(\overline{{\mathbb S}^{N}_+})$ with respect to the norm
\begin{equation}\label{eq:h1sfera}
\|\psi\|_{H^1({\mathbb S}^{N}_+;\theta_1^{1-2s})}=\bigg(
\int_{{\mathbb S}^{N}_+}\theta_1^{1-2s}\big(|\nabla_{{\mathbb
        S}^{N}}\psi(\theta)|^2+\psi^2(\theta)\big)dS
\bigg)^{\!\!1/2}.
\end{equation}
For every $a\in L^\infty(\S^{N-1})$, let
\begin{equation}\label{firsteig}
  \mu_1(a):=\min_{\psi\in H^1({\mathbb
  S}^{N}_+;\theta_1^{1-2s}) \setminus\{0\}}\frac{ \int_{{\mathbb S}^{N}_+}\theta_1^{1-2s}|\nabla
\psi |^2\,dS-
    \kappa_s \int_{{\mathbb
      S}^{N-1}}a  \psi^2 \,dS'}{\int_{{\mathbb
    S}^{N}_+}\theta_1^{1-2s}\psi^2 \,dS},
\end{equation}
where
\begin{equation*}
\kappa_s=\frac{\Gamma(1-s)}{2^{2s-1}\Gamma(s)}.
\end{equation*}
The quantity $\mu_1(a)$ is an eigenvalue appearing from a
change of  polar coordinates in some Dirichlet energy defined on the
half-space $\RNp$, see \cite{FF1}. In \cite{FF1}  it is also observed
that the operator
$\Dsm-a(x/|x|) |x|^{-2s}$ is positive definite provided
 \be \label{eq:positivcond}
\mu_1(a)+\frac{(N-2s)^2}{4}>0,\ee
 see also Lemma
\ref{l:hardyboundary} below. Throughout   this
paper, we will always assume
\eqref{eq:positivcond}.

 The following theorem gives conditions on the coefficient $a$ for
 essential self-adjointness of the operator $ \Dsm-a(x/|x|)
 |x|^{-2s}$.

\begin{Theorem}\label{t:esa}
Assume that $s\in(0,1)$, $m\geq0$, $N>2s$, and $a\in L^\infty(\S^{N-1})$ with
\[
\mu_1(a)+\frac{(N-2s)^2}{4}>0.
\]
 Then the  operator
 \[
\Dsm-a(x/|x|)
|x|^{-2s}\quad\text{with domain $C^\infty_c(\R^N\setminus\{0\})$}
\]
 is
essentially self-adjoint in $L^2(\R^N)$ if and only if
\begin{equation}\label{eq:11}
-\mu_1(a)\leq
\frac{(N-2s)^2}{4}-s^2.
\end{equation}
\end{Theorem}

\noindent  If
  $a$ is constant then $\mu_1(a)$ can be obtained implicitly from the usual Gamma function. Indeed, pick   $\a\in\left(0,
\frac{N-2s}{2}\right)$  and let
\begin{equation}\label{eq:lamal}
\l(\a)=2^{2s}  \frac{\Gamma\left(\frac{N+2s+2\a}{4}\right)}{\Gamma\left(\frac{N-2s-2\a}{4}\right)}
\frac{\Gamma\left(\frac{N+2s-2\a}{4}\right)}{\Gamma\left(\frac{N-2s+2\a}{4}\right)}.
\end{equation}
From \cite[Proposition 2.3]{FF}, we have
$$
\mu_1(\l(\a))=\a^ 2-\left( \frac{N-2s}{2} \right)^{\!\!2}\quad \text{for
  all }\a\in   \left(0, \frac{N-2s}{2}\right).
$$
Therefore,  by combining Theorem
\ref{t:esa} and \cite[Proposition 2.3]{FF}, we obtain the following corollary.
\begin{Corollary}
Let  $\a\in\left(0, \frac{N-2s}{2}\right]$. Then
$\Dsm-\l(\a) |x|^{-2s}$ with domain
$C^\infty_c(\R^N\setminus\{0\})$ is essentially self-adjoint in
$L^2(\R^N)$  if and only if $  \a\geq s$.
\end{Corollary}
\noindent Since the map $\a\mapsto \l(\a)$ is decreasing, we also obtain the
following corollary.
\begin{Corollary}\label{cor:th-ess-l} Let $\b\in\R$ and $\lambda$ be given by \eqref{eq:lamal}. Then
the operator $\Dsm-\b |x|^{-2s}$ with domain
$C^\infty_c(\R^N\setminus\{0\})$ is essentially self-adjoint in
$L^2(\R^N)$  if and only if $  \b\leq  \l( s)$.
\end{Corollary}
Let us note that
$$
\l(s)=2^{2s}\frac{ \G(\frac{N+4s}{4})}{\G(\frac{N-4s}{4})}.
$$
If $s=1/2$ and $N=3$ then $\l(1/2)=1/2$. In this case, the
essential self-adjointness below the threshold $\l(1/2)=1/2$ was known, see Kato \cite{Kato} and
Herbst \cite{Herbst}; moreover, the sharpness of the threshold $1/2$  in
3 dimensions was obtained in \cite[Corollary 1]{LOR}.
 In higher dimensions $N\geq 3$ and for  $s=1/2$,
Ichinose in \cite{Ich} proved  essential self-adjointness of
 $\Dsm-\b |x|^{-2s}$
provided
$\b < \frac{N-2}{2}$ using the Kato-Rellich perturbation result and
the Hardy inequality.  Our result in Corollary \ref{cor:th-ess-l}
improves the results in \cite{Ich} because, for $s=1/2$, we have $\l(1/2)=
\frac{N-2}{2}$. In addition, we also obtain the sharpness of the
threshold $\l(1/2)=\frac{N-2}{2}$, thus extending \cite[Corollary
1]{LOR} to higher dimensions.

We remark that the precise threshold for non essential
self-adjointness in the local case $s=1$  is
$\l(1)=\frac{(N-2)^2}{4}-1$; we refer to \cite{FMT-jfa,KSWW,Simon}
for such local case. We also mention that the case of $a$ not
constant was treated in the non-relativistic case in \cite{FMT}.

We observe that as a direct consequence of Theorem \ref{t:esa} and
Kato-Rellich Perturbation Theorem the operator
 $\Dsm-m^{2s}-a(x/|x|)
|x|^{-2s}$ with domain $C^\infty_c(\R^N\setminus\{0\})$ is
essentially self-adjoint in $L^2(\R^N)$ if and only if $  -\mu_1(
a)\leq \frac{(N-2s)^2}{4}-s^2$. We refer to \cite{CMS} for the study
of asymptotics of the eigenstates of relativistic operators  of type $\Dsm-m^{2s}$.

\bigskip Our argument for proving essential self-adjointness is quite direct
and it is inspired by \cite{Simon}. The proof is based on
a contradiction argument as follows: if $A=\Dsm-\frac{a(x/|x|)}{|x|^{2s}}$, with domain $C^\infty_c(\R^N\setminus\{0\})$, is not
essentially self adjoint in $L^2(\R^N)$ then $\Ds-\frac{a(x/|x|)}{|x|^{2s}}$, with domain $C^\infty_c(\R^N\setminus\{0\})$, is
not essentially self adjoint in $L^2(\R^N)$ as well (by Kato-Rellich),
see Lemma \ref{eq:meq0}.

 Therefore there exists a function $w\in
L^2(\R^N)$, $w\neq0$, such that $\Ds w-a(x/|x|) |x|^{-2s}w+w=0 $ in the
sense of distributions, namely
\begin{equation}\label{eq:essttp}
\int_{\R^N}[\Ds \vp -a(x/|x|)
|x|^{-2s}\vp +\vp]w =0, \quad \textrm{for all } \vp \in C^\infty_c(\R^N\setminus\{0\}).
\end{equation}
The idea is now to construct appropriate test functions in \eqref{eq:essttp} to get $w\equiv 0$ leading to  a contradiction. In order to do this, we will approximate $a$ by smooth functions $a_n\in C^\infty(\S^{N-1})$
such that $\mu_1(a_n)\to \mu_1(a)$.

The test functions we will consider are then solutions to the  partial differential equations
$$
\Ds v_n + v_n -a_n(x/|x|)
|x|^{-2s}v_n= f,
$$
for arbitrary $f\in C^\infty_c(\R^N\setminus\{0\})$, $f\geq0$,
$f\neq0$. Such functions $v_n$ enjoy the following estimates at the
origin and at infinity:
$$
v_n\leq C |x|^{\g_n} \textrm{ in } B_{r_0},\qquad v_n\leq C |x|^{\a_n} \textrm{ in } \R^N\setminus B_{r_0},
$$
where $B_{r_0}=\{x\in\R^N:|x|<r_0\}$, $r_0>0$, and
$$
\g_n=-\frac{N-2s}{2}+\sqrt{\bigg(\frac{N-2s}
    {2}\bigg)^{\!\!2}+\mu_{1}(a_n)},   \quad     \a_n= -\frac{N-2s}{2}-\sqrt{\bigg(\frac{N-2s}
    {2}\bigg)^{\!\!2}+\mu_{1}(a_n)}.
$$
Regularity theory implies that $v_n\in
C^\infty(\R^N\setminus\{0\})$.

 We then cut  off the $v_n$'s at the origin
and at infinity, in order to use them as test functions in
\eqref{eq:essttp}. Then, thanks to \eqref{eq:11} and some tricky
integration by parts in the nonlocal framework, we end up with
$$
\int_{\R^N}wf\, dx=\int_{\R^N}(a(x/|x|)-a_n(x/|x|) )|x|^{-2s}  v_nw\,dx.
$$
Finally, from the above estimates of $v_n$ it follows that
$\int_{\R^N}wf\, dx=0 $, thus contradicting that $w\neq 0$. This
program is elaborated in details in Section
\ref{sec:essent-self-adjo}.

We observe that the above described  arguments can be adapted to treat operators
of the type  $\Dsm-a(x/|x|)
 |x|^{-2s}+h(x)$ where $h\in L^\infty_{loc}(\R^N)\cap L^p(B_r)$, for
 some $p >N/(2s)$, $r>0$, and $h$  is bounded in a
neighborhood of $\infty$, see Remark \ref{rem:self-h}.

\bigskip
To prove non essential self-adjointness of a densely defined operator,
it is generally inevitable to solve some partial differential
equations (mostly, boundary value eigenvalue problems) for which the
solutions are known explicitly or at least have some qualitative
properties that can be handled.  In our situation, we would like to
prove that $\Dsm -a(x/|x|) |x|^{-2s}$ is not essentially self-adjoint
when $ -\mu_1(a)> (N-2s)^2/{4}-s^2$.

 To show this we argue by
contradiction and assume that $\Dsm -a(x/|x|) |x|^{-2s}$ is
essentially self-adjoint, which is equivalent to the density of the
range of  $\Dsm -a(x/|x|) |x|^{-2s}+d$ in $L^2(\R^N)$ for all $d>0$.

We show in Lemma
\ref{eq:meq0} that this is equivalent to $(-\Delta+ m^2+
b)^s-a(x/|x|) |x|^{-2s}$ having dense range in $L^2(\R^N)$ for all
$b>0$;
this circumstance is  ruled out by constructing a function $f\neq 0$ solving
the equation $(-\Delta+ m^2+ b)^s f-a(x/|x|) |x|^{-2s}f=0 $ in $\R^N$
with $f\in L^2(\R^N)$ provided $ -\mu_1(a)>
\frac{(N-2s)^2}{4}-s^2$. The advantage of considering $(-\Delta+ m^2+
b)^s $ instead of $ \Dsm +d$ is the exponential decay at infinity of the fundamental solution of the
former operator, which is crucial in
our analysis and which fails for the latter operator for $m=0$.
This argument  will be developed in details in Section \ref{s:non-ess}.

\section{Some preliminaries and Notations}\label{s:sp}
We start by recalling the integral representation of $\Dsm$: for every $u\in C^2_c(\R^N)$
$$
\Dsm u(x)={c_{N,s}}m^{\frac{N+2s}{2} }
P.V.\int_{\R^{N}}\frac{u(x)-u(y)}{|x-y|^{\frac{N+2s}{2}} } K_{ \frac{N+2s}{2}}(m|x-y|)\,dy+m ^{2s} u(x),
$$
where $m\geq0$ and
$$
c_{N,s}=
{2^{-(N+2s)/2+1}}\pi^{-\frac
  N2}2^{2s}\frac{s(1-s)}{\Gamma(2-s)},
$$
see \cite{FF1}. The kernel $K_{\nu}$  denotes the modified Bessel function of the second
kind with order $\nu$.
We recall that, for $\nu>0$,
\be\label{eq:decKnu0}
 K_{\nu}(r)\sim
\frac{\G(\nu)}{2} \left(\frac{r}{2} \right)^{-\nu}
\ee
 as $r\to 0$
and $K_{-\nu}=K_{\nu}$ for $\nu<0$, while
\be\label{eq:decKnuInf}
K_{\nu}(r)\sim \frac{\sqrt{\pi }}{\sqrt{2}} r^{-1/2}e^{-r}\ee as
$r\to +\infty$, see
\cite{El}. Furthermore there holds
 $$
K_\nu'(r)=-\frac{\nu}{r}K_\nu(r) -K_{\nu-1}(r).
$$
 The Dirichlet form associated to
$\Dsm$ on $C^\infty_c(\R^N)$ is given by
\begin{align}\label{eq:3}
  (u,v)_{
    H^{s}_m(\R^N)}:&=\int_{\R^N}(|\xi|^2+m^2)^s\widehat{u}(\xi)\ov{\widehat{v}(\xi)} d\xi\\
    &\nonumber =\frac{c_{N,s}}{2}m^{\frac{N+2s}{2} }
\int_{\R^{2N}}\frac{(u(x)-u(y))(v(x)-v(y))}{|x-y|^{\frac{N+2s}{2}} } K_{ \frac{N+2s}{2}}(m|x-y|)\,dx\,dy\\
\nonumber & \hspace{1cm} +m ^{2s} \int_{\R^N}u(x)v(x)dx,
\end{align}
where $\widehat{u}$ denotes the unitary Fourier transform of $u$.
We  define $ H^{s}_m(\R^N)$ as the completion of
$C^\infty_c(\R^N)$ with respect to the norm induced by the
scalar product \eqref{eq:3}. If $m>0$, $ H^{s}_m(\R^N)$ is nothing
but the standard $ H^{s}(\R^N)$; then, we will write  $H^{s}(\R^N)$
without the subscript ``$m$''.\\

The operator $\Dsm$ enjoys an extension property reminiscent of
the Caffarelli-Silvestre extension \cite{CSilv}, see \cite{FF1}.  Let us recall it
via the Bessel Kernel which is given by
\be\label{eq:Bess-Kernel}
{P}_m(z)=C_{N,s}'\, t^{2s}m^{\frac{N+2s}{2} }
|z|^{-\frac{N+2s}{2}}K_{\frac{N+2s}{2}}(m|z|),
\ee
with $z=(t,x)\in\R\times\R^N$ and  some normalization constant  $C'_{N,s}$. Pick   $u\in H^s_m(\R^N)$ and set
 $$
 w(t,x)=(P_m(t,\cdot)*u)(x).
 $$
Then, see \cite{FF1}, we have  that $w\in H^1_m(\RNp;t^{1-2s})$ and moreover
 \begin{equation}\label{eq:ext}
\begin{cases}
- \div(t^{1-2s}\n w)(t,x)+m^2t^{1-2s} w(t,x)=0,&\textrm{in  } \RNp,\\
 -\lim_{t\to 0} t^{1-2s}\frac{\de w}{\de t}(t,x)
 = \k_s\Dsm u (x),&\textrm{on  } \R^N,
\end{cases}
\end{equation}
in a weak sense, where $\R^{N+1}_+ =\{z=(t,x):t\in(0,+\infty),\
x\in\R^N\}$. Here $H^1_m(\RNp;t^{1-2s})$ is the completion of
$C^\infty_c(\ov{\RNp})$ with respect to the norm
$\int_{\RNp}t^{1-2s}|\n w|^2\,dtdx+m^2\int_{\R^N}t^{1-2s}w^2\,dtdx$.
In fact, performing Fourier transform in the above equations, we can
see that the Bessel Kernel $P_m(t,x)$ is the Fourier transform of
$\xi\mapsto \vartheta(\sqrt{|\xi|^2+m^2}t)$, where $\vartheta(r)=\frac{2}{\G(s)}\left( \frac{r}{2}\right)^s\,K_s(r)$ solves
$$
\begin{cases}
\vartheta''+\frac{(1-2s)}{t}\vartheta'-\vartheta=0,\\
\vartheta(0)=1.
\end{cases}
$$
This then implies that
\be\label{eq:intPm} \int_{\R^N}P_m(t,x)dx=\vartheta(mt). \ee

Due to homogeneity properties of problem
\eqref{eq:ext}, we are naturally lead to consider an
angular eigenvalue problem. Let $H^1({\mathbb
  S}^{N}_+;\theta_1^{1-2s})$ be defined as in \eqref{eq:h1sfera}.
Since the weight $\th^{1-2s}_1$ belongs to the second Muckenhoupt 
class, the embedding
\[
H^1({\mathbb
  S}^{N}_+;\theta_1^{1-2s})\hookrightarrow \hookrightarrow
L^2({\mathbb S}^{N}_+;\theta_1^{1-2s})
\]
 is compact, where
\[
L^2({\mathbb
  S}^{N}_+;\theta_1^{1-2s}):=\Big\{\psi:{\mathbb S}_+^{N}\to\R\text{
  measurable such that }{\textstyle{\int}}_{{\mathbb
    S}^{N}_+}\theta_1^{1-2s}\psi^2(\theta)\,dS<+\infty\Big\}.
\]
Letting   $a\in L^q(\S^{N-1})$, for some $q>N/(2s)$,  the first
eigenvalue of the angular component of the extended operator
 $$
   \mu_1(a)=\min_{\psi\in H^1({\mathbb
  S}^{N}_+;\theta_1^{1-2s}) \setminus\{0\}}\frac{ \int_{{\mathbb S}^{N}_+}\theta_1^{1-2s}|\nabla
\psi |^2\,dS-
    \kappa_s \int_{{\mathbb
      S}^{N-1}}a  \psi^2 \,dS'}{\int_{{\mathbb
    S}^{N}_+}\theta_1^{1-2s}\psi^2 \,dS}
$$
is attained by an eigenfunctions $\psi$  which  does not change sign and satisfies
\begin{equation}\label{eq:1eig}
 \begin{cases}
    -\dive\nolimits_{{\mathbb S}^{N}}(\theta_1^{1-2s}\nabla_{{\mathbb
        S}^{N}}\psi)=\mu_1(a)\,
    \theta_1^{1-2s}\psi, &\text{in }{\mathbb S}^{N}_+,\\[5pt]
-\lim_{\theta_1\to 0^+} \theta_1^{1-2s}\nabla_{{\mathbb
    S}^{N}}\psi\cdot {\mathbf
  e}_1=\kappa_s a(\theta')\psi,&\text{on }\partial {\mathbb S}^{N}_+={\mathbb S}^{N-1}.
  \end{cases}
\end{equation}
The following result is essentially  contained in \cite{FF1}.
\begin{Lemma} \label{l:hardyboundary}
Let $q>N/(2s)$ and $a\in L^q(\S^{N-1})$ such that
$$
\mu_1(a)+ \bigg(\frac{N-2s}2\bigg)^{\!\!2}>0.
$$
 Then there exists a constant $C_{a,N,s}>0$ such that, for all $w\in H^1_0(\RNp;t^{1-2s})$,
 $$
  \int_{\RNp} t^{1-2s}|\nabla w|^2\,dt\,dx- \kappa_s \int_{\R^N}
    \frac{a(x/|x|)}{|x|^{2s}}w^2\,dx  \\
    \geq C_{a,N,s}
    \int_{\RNp} t^{1-2s}|\nabla w|^2\,dt\,dx.
    $$
Equivalently, we have that
$$
    \int_{\R^N} |\xi| ^{2s}\widehat{ \vp}^2\,d\xi-   \int_{\R^N}
    \frac{a(x/|x|)}{|x|^{2s}}\vp^2\,dx
    \geq C_{a,N,s}
     \int_{\R^N} |\xi| ^{2s}\widehat{ \vp}^2\,d\xi
  $$
    for all   $\vp\in H^s_0(\R^N)$.
\end{Lemma}
\begin{Remark}\label{rem:CaNscont}
It is useful to remark that the best constant $C_{a,N,s}$ in Lemma
\ref{l:hardyboundary}  depends
continuously on $a$ as a  mapping in $ L^q(\S^{N-1})$, see Remark 2.5
in \cite{FF1}.
\end{Remark}

\noindent We will also need the following result from \cite{FW}.
\begin{Lemma}[\cite{FW}, Lemma 2.1] \label{lem:FW}
Let $\O$ be a bounded open set.  Then there
 exists a positive constant $C=C(N,s,\O)>0$ such that for all  $\vp\in C^2_c(\O)$ and for all  $x\in \R^N$
$$
\left|\Ds \vp (x) \right|\leq \frac{C\|\vp\|_{C^2(\R^N)}}{1+|x|^{N+2s}} .
$$

\end{Lemma}

\section{Essentially self-adjointness}\label{sec:essent-self-adjo}

In this section we shall prove that the operator $A'=\Dsm-a(x/|x|)
|x|^{-2s}$ with domain $C^\infty_c(\R^N\setminus\{0\})$ is
essentially self-adjoint in $L^2(\R^N)$ provided
\[
 -\mu_1(a)\leq
\frac{(N-2s)^2}{4}-s^2.
\]
 This stands to be a generalization of the
case $s=1$ by Kalf, Schmincke, Walter,
W\"{u}st \cite {KSWW}, see also Simon \cite{Simon},  and of the case
in which $s=1/2$ and $a$ is   constant which was treated by Kato
\cite{Kato}, see also Herbst \cite{Herbst}.

 For the proof, we will need some
technical lemmata. Let us first observe that it is not restrictive to take $m=0$.
  \begin{Lemma}\label{eq:meq0}
For $s\in(0,1)$, $V\in L^2_{loc}(\R^N\setminus\{0\})$ and $b>0$, let
us consider  $A=\Ds-V$ and $B=(-\D +b)^s-V$ with domain
$C^\infty_c(\R^N\setminus\{0\})$.
\begin{enumerate}[$(i)$]
\item $A$ is essentially self-adjoint
on $L^2(\R^N)$ if and only if $B$ is essentially self-adjoint on
$L^2(\R^N)$.
\item If there exists  $C>0$ such that
\begin{equation}\label{eq:pot}
\int_{\R^N}(|\xi|^2+b)^{s}\widehat{\vp}^2(\xi)\,d\xi-\int_{\R^N} V(x)
\vp^2(x)\,dx\geq C \|\vp\|^2_{H^s(\R^N)} \quad
\text{for every } \vp\in C^\infty_c(\R^N\setminus\{0\}),
\end{equation}
then $B$ is essentially self-adjoint on
$L^2(\R^N)$ if and only if $B$ has dense  range in $L^2(\R^N)$.
\end{enumerate}
\end{Lemma}
 \proof
To prove $(i)$, we observe that, by
Fourier transform
and Parseval identity,
 $$
 \|(B-A)u\|^2_{L^2(\R^N)}=\int_{\R^N}\big[(|\xi|^2+b)^s-|\xi|^{2s}
 \big]^2|\widehat{u}(\xi)|^2d\xi,
 $$
for all $u\in C^\infty_c(\R^N\setminus\{0\})$.
Using
the elementary inequality $0\leq (a+b)^s-a^s\leq
b^s$, which holds for every $a,b\in[0,+\infty)$ and $s\in(0,1)$, it follows that
 $$
 \|(B-A)u\|^2_{L^2(\R^N)}\leq b^{2s}\int_{\R^N}  |\widehat{u}(\xi)|^2d\xi=b^{2s}\|u\|^2_{L^2(\R^N)}.
 $$
 Therefore,
for $q\in (0,1)$ we get
  \begin{align*}
 & \|(B-A)u\|_{L^2(\R^N)}\leq  q \|A u\|_{L^2(\R^N)}+
  b^{s}\|u\|_{L^2(\R^N)}, \\
&  \|(B-A)u\|_{L^2(\R^N)}\leq  q \|B u\|_{L^2(\R^N)}+
  b^{s}\|u\|_{L^2(\R^N)},
\end{align*}
for all $u\in C^\infty_c(\R^N\setminus\{0\})$, i.e. $B-A$ is both
$A$-bounded and $B$-bounded with relative bound $q<1$. Then by the Kato-Rellich Theorem (see e.g. \cite[Theorem
 X.12]{reedsimon}) it follows that if $A$ is essentially self-adjoint
 then $B=A+(B-A)$ is
 essentially self-adjoint; in the same way, if $B$ is essentially self-adjoint
 then $A=B+(A-B)$ is
 essentially self-adjoint, thus proving $(i)$.

We recall (see e.g. \cite[Theorem X.26]{reedsimon}) that if a symmetric operator  is  strictly
positive, then it is essentially self-adjoint if and only if its range
is dense.
 Since by assumption \eqref{eq:pot} $B$  is a strictly positive
symmetric operator, we deduce statement $(ii)$.
\QED

\begin{Remark}\label{sec:essent-self-adjo-1}
Let us observe that our potential $V(x)=\frac{a(x/|x|)}{|x|^{2s}}$
satisfies \eqref{eq:pot} for every $b>0$ provided condition \eqref{eq:positivcond} is satisfied.
Indeed,  for $\e>0$ and $u\in H^s(\R^N)$, we have
\begin{align*}
  \int_{\R^N}& (|\xi|^2+b)^s
  |\widehat{u}|^2(\xi)d\xi-\int_{\R^N}\frac{a(x/|x|)}{|x|^{2s}}u^2dx\\
  &= (1-\e) \int_{\R^N} (|\xi|^2+b)^s |\widehat{u}|^2(\xi)d\xi
  +\e  \int_{\R^N} (|\xi|^2+b)^s |\widehat{u}|^2(\xi)d\xi -\int_{\R^N}\frac{a(x/|x|)}{|x|^{2s}}u^2dx\\
  &= (1-\e)\bigg[ \int_{\R^N} (|\xi|^2+b)^s |\widehat{u}|^2 d\xi-
  \int_{\R^N} \frac{a_\e(x/|x|)}{|x|^{2s}} u^2dx \bigg] +\e
  \int_{\R^N} (|\xi|^2+b)^s |\widehat{u}|^2 d\xi,
\end{align*}
where $a_\e=\frac{a}{1-\e}$. By continuous dependence of $\mu_1$ on
$a$ and \eqref{eq:positivcond}, there exists $\e_0=\e_0(a,N,s)>0$ such
that
 $$
 \mu_1(a_\e)+\frac{(N-2s)^2}{4}>0
 $$
 for all $\e\in (0,\e_0)$.
 By Lemma \ref{l:hardyboundary}, we have
 $$
    \int_{\R^N}  (|\xi|^2+b)^s\widehat{ u}^2\,d\xi-   \int_{\R^N}
    \frac{a_\e(x/|x|)}{|x|^{2s}}u^2\,dx
    \geq C_{a_\e,N,s}
     \int_{\R^N} |\xi| ^{2s}\widehat{ u}^2\,d\xi.
  $$
Therefore
\begin{align*}
 \int_{\R^N} (|\xi|^2+b)^s &|\widehat{u}|^2
 d\xi-\int_{\R^N}\frac{a(x/|x|)}{|x|^{2s}}u^2dx\\
&\geq(1-\e) C_{a_\e,N,s}
     \int_{\R^N} |\xi| ^{2s}\widehat{ u}^2\,d\xi+\e  \int_{\R^N} (|\xi|^2+b)^s |\widehat{u}|^2 d\xi.
\end{align*}
Hence,  by Parseval identity, for all $\e\in (0,\e_0)$ we have that
\begin{align*}
 \int_{\R^N} (|\xi|^2+b)^s |\widehat{u}|^2 d\xi&-\int_{\R^N}\frac{a(x/|x|)}{|x|^{2s}}u^2dx\geq(1-\e) C_{a_\e,N,s}
     \int_{\R^N} |\xi| ^{2s}\widehat{ u}^2\,d\xi+\e b^s  \int_{\R^N}    {u}^2 dx.
\end{align*}
\end{Remark}

\medskip
\noindent The following uniforms decay estimates will be useful in the sequel.
\begin{Lemma}\label{l:stima}
Let
$a_n\in C^\infty(\S^{N-1})$ be such that $a_n\to a$ in $L^q(\S^{N-1})$, for some $q>N/(2s)$, and
$\mu_1(a_n)\to \mu_1(a)$ as $n\to\infty$. Assume that
$$
\mu_1(a)+ \bigg(\frac{N-2s}2\bigg)^{\!\!2}>0.
$$
 Let  $v_n\in H^s_0(\R^N)$ be a sequence of  functions such that $v_n>0$
a.e. in $\R^N$ and 
$\{v_n\}_n$ is bounded in $ {H^s_0}(\R^N)$.
\begin{enumerate}[(i)]
\item If 
 \begin{equation*}
  \Ds v_n -a_n(x/|x|) |x|^{-2s}v_n \leq 0,\quad\text{in }B_R,
\end{equation*}
for some $R>0$, then there exist $C>0$ and $r_0\in(0, R)$ (independent of $n$) such that
\be\label{eq:vxest} 
v_n(x)\leq C|x|^{\gamma_n} \textrm{ for a.e. $x\in {B_{r_0}}$}
\ee
where
\[
\gamma_n:=-\frac{N-2s}{2}+\sqrt{\bigg(\frac{N-2s}
    {2}\bigg)^{\!\!2}+\mu_{1}(a_n)}.
\]
\item  If 
 \begin{equation*}
  \Ds v_n -a_n(x/|x|) |x|^{-2s}v_n \leq 0,\quad\text{in }\R^N\setminus
  B_R,
\end{equation*}
for some $R>0$, then there exist $C>0$ and $r_0> R$ (independent of $n$) such that
\be\label{eq:vxestinfty}
v_n(x)\leq C|x|^{\a_n} \textrm{ for a.e. $x\in \R^N\setminus {B_{r_0}}$}
\ee
  for $n$ sufficiently large,  where 
$$
\a_n= -\frac{N-2s}{2}-\sqrt{\bigg(\frac{N-2s}
    {2}\bigg)^{\!\!2}+\mu_{1}(a_n)}.
$$
\end{enumerate}
\end{Lemma}
\proof
To prove $(i)$,  let $w_n\in  H^1_0(\R^{N+1}_+;t^{1-2s})$ be the
Caffarelli-Silvestre extension of $v_n$, so that $w_n$ solves
\begin{equation}\label{eq:5}
\begin{cases}
\dive(t^{1-2s}\n w_n)=0,&\textrm{ in }\R^{N+1}_+,\\
-  \lim \limits_{t \to 0^+}t^{1-2s}\,\frac{\partial w_n}{\partial t}\leq
\k_s a_n(x/|x|) |x|^{-2s}v_n , &\textrm{ on } B_R.
\end{cases}
\end{equation}
If $r_0\in(0,R)$, from the regularity estimates in \cite{JLX} (see
also \cite[Proposition 3.3]{FF1}),  we deduce that
$w_n|_{S_{r_0}^+}$ is uniformly bounded, where $S_{r_0}^+=\{z\in
\R^{N+1}_+:|z|=r_0\}$.

 Hence there exists $C>0$
independent of $n$ such that $0\leq w_n\leq C\tilde w_n$ on
$S_{r_0}^+$ for $n$ sufficiently large, where $\tilde
w_n(z):=|z|^{\gamma_n}\psi_n(z/|z|)$ with $\psi_n$ being the positive
$L^2({\mathbb S}^{N}_+;\theta_1^{1-2s})$-normalized
eigenfunction corresponding
to $\mu_1(a_n)$. Since $\tilde w_n$ solves
\begin{equation}\label{eq:6}
\begin{cases}
\dive(t^{1-2s}\n \tilde w_n)=0,&\textrm{ in }\R^{N+1}_+,\\
-  \lim \limits_{t \to 0^+}t^{1-2s}\,\frac{\partial \tilde w_n}{\partial t}=
\k_s a_n(x/|x|) |x|^{-2s}\tilde w_n, &\textrm{ on }\R^N,
\end{cases}
\end{equation}
testing the difference between \eqref{eq:5} and \eqref{eq:6}
multiplied by $C$ with $(w_n-C\tilde w_n)^+$, integrating by parts,
using that $v_n>0$,  and
invoking Lemma \ref{l:hardyboundary}, we obtain that $w_n\leq C\tilde
w_n$  a.e. in $B_{r_0}^+=\{z\in
\R^{N+1}_+:|z|<r_0\}$. Hence $v_n(x)\leq C
|x|^{\gamma_n}\psi_n(x/|x|)$  for a.e. $x\in B_{r_0}$ and the conclusion
follows from an uniform upper bound of $\psi_n$ (which follows e.g.
from \cite[Proposition 3.3]{FF1}).\\

To prove $(ii)$, we consider the  Kelvin transform of $v_n$ given by
$\ti{v}_n=|x|^{2s-N}v_n(x/|x|^2)$. We have that $ \ti{v}_n\in {H^s_0}(\R^N)$
with $ \|\ti{v}_n\|_{ {H^s_0}(\R^N)}=\|{v}_n\|_{{H^s_0}(\R^N)} $
(see \cite[Lemma 2.2]{FW}) and
$$
 \Ds \ti{v}_n -a_n(x/|x|) |x|^{-2s}\ti{v}_n \leq 0\quad\text{in }B_{1/R}.
$$
From $(i)$,  
 for some $C_1>0$ and $r_0>R$ (independent on $n $), we have 
$$
0\leq \ti{v}_n(x) \leq C_1 |x|^{\g_n},  \quad
\text{for all } x\in
B_{r_0}\setminus\{0\}
$$
which yields
$$
v_n(x)\leq C |x|^{\alpha_n},   \quad
\text{for all }  x\in\R^N\setminus
B_{r_0},
$$
where $\a_n= -\frac{N-2s}{2}-\sqrt{\big(\frac{N-2s}
    {2}\big)^{2}+\mu_{1}(a_n)}$.
\QED

\begin{Theorem}\label{th:essself}
Assume that $s\in (0,1)$, $N>2s$, and $a\in L^q({\mathbb
      S}^{N-1})$ for some   $q> \max\big(\frac{N}{2s},2\big)$. Then the operator
 $A=\Ds-a(x/|x|)
|x|^{-2s}$ with domain $C^\infty_c(\R^N\setminus\{0\})$ is
essentially self-adjoint in $L^2(\R^N)$ provided $  -\mu_1( a)\leq
\frac{(N-2s)^2}{4}-s^2.$
\end{Theorem}

 \proof
 The proof of the theorem will be separated into two cases.\\
\noindent
 \textbf{Case 1: } \be\label{eq:musHmeq}
 -\mu_1(a)<
\frac{(N-2s)^2}{4}-s^2.
  \ee
 By Lemma \ref{l:hardyboundary}  we have that
$$
( A \vp,\vp)_{L^2(\R^N)}=\int_{\R^N}(A\varphi)(x)\varphi(x)\,dx\geq 0,
\quad\text{for all } \vp\in
C^\infty_c(\R^N\setminus\{0\}),
$$
so that $A$  is nonnegative definite in
$C^\infty_c(\R^N\setminus\{0\})$.

From the Kato-Rellich Theorem and well-known self-adjointness criteria for positive operators (see
\cite[Theorem X.26]{reedsimon}), $A$ with domain $C^\infty_c(\R^N\setminus\{0\})$ is
essentially self-adjoint in $L^2(\R^N)$ if and only if $\mathop{\rm
  Range}(A+1)=(A+1)(C^\infty_c(\R^N\setminus\{0\}) )\subset L^2(\R^N)$ is
dense in $L^2(\R^N)$.

We argue by contradiction and
assume that $A$ is not essentially self-adjoint so that
\[
(A+1)(C^\infty_c(\R^N\setminus\{0\}) )\subset L^2(\R^N)
\]
 is not
dense in $L^2(\R^N)$. Then there exists $w\in L^2(\R^N)$,
$w\neq0$ such that $(w, u)_{L^2(\R^N)}=0 $ for all $u \in
(A+1)(C^\infty_c(\R^N\setminus\{0\}) )$. In particular
\begin{equation}\label{eq:7}
  \int_{\R^N} w(x)\bigg[\Ds \vp(x)-a\big({x}/{|x|}){|x|^{-2s}} \vp (x)+\vp(x)\bigg]\,dx=0\quad \text{for all }\vp\in
  C^\infty_c(\R^N\setminus\{0\}),
\end{equation}
i.e.
$$
  \Ds w-a(x/|x|) |x|^{-2s}w +w =0 \qquad \textrm{ in } \calD'(\R^N\setminus\{0\}).
$$
 We will reach a contradiction by showing that $w\equiv 0$; to prove
 that $w\equiv 0$, we will prove that $\int_{\R^N}fw dx =0$ for every $f\in C^\infty_c(\R^N\setminus\{0\})$, $f\geq0$,
$f\neq0$. 
To this aim, let us fix 
\[
f\in C^\infty_c(\R^N\setminus\{0
\})\quad\text{such that}\quad 
 f\geq 0,\quad f\not\equiv0.
\]

\medskip\noindent
{\bf Step 1.} 
By density, there exists a sequence  $a_n\in C^\infty(\S^{N-1})$ such that
 \be\label{eq:antoa}
 a_n\to a \quad  \textrm{in } L^{q}(\S^{N-1})\quad  \textrm{as } n\to \infty.
 \ee
 Then by \cite[Lemma 2.1]{FF1}, we have that   $\mu_{1}(a_n)\to \mu_{1}(a)$.
By  \eqref{eq:musHmeq}, we have that, for some $\e_0>0$,
\be\label{eq:muanHarn}
-s+\sqrt{\bigg(\frac{N-2s}
    {2}\bigg)^{\!\!2}+\mu_{1}(a_n)}>\e_0>0
\ee
for every large $n$.

By the Lax-Milgram
theorem, for every $n\in\N$, there exists $v_n\in H^s(\R^N)$ such that
 \be\label{eq:vstf}
  \Ds v_n -a_n(x/|x|) |x|^{-2s}v_n+v_n=f .
\ee
Multiplying \eqref{eq:vstf} by the negative part of $v_n$ and using Lemma
\ref{l:hardyboundary}, we can see that $v_n\geq 0$ since $f\geq 0$.
By the Harnack inequality  $v_n >0$ in $\R^N\setminus\{0\}$  (see \cite{JLX}). By \eqref{eq:muanHarn} and Remark \ref{rem:CaNscont},  the sequence  $(v_n)_n$
is bounded in $H^s(\R^N)$.  By regularity theory
$v_n\in C^\infty(\R^N\setminus\{0\})$ (see also \cite{JLX}). Moreover   Lemma
\ref{l:stima},  and \eqref{eq:muanHarn} imply  that  $\Ds v_n\in L^2(\R^N)$ so that $v_n\in H^{2s}(\R^N)$.

\medskip\noindent
{\bf Step 2.}  Let $\eta\in C^\infty_c(\R)$ be such that $0\leq \eta\leq1$, $\eta(t)=1$ for  $|t|\leq 1$ and
$\eta(t)=0$ for $|t|\geq 2$. We put $\eta_\d(t)=\eta(\frac{t}{\d})$
and $\eta_R(t)=\eta(\frac{t}{R})$ so that $(1-\eta_\d)  \eta_R
v_n\in C^\infty_c(\R^N\setminus\{0\})$. We put
$v_{n,\d}=(1-\eta_\d)v_n $. Then from \eqref{eq:7} we have that
\begin{equation}\label{eq:vndR}
  \int_{\R^N}w\Big(\Ds( \eta_R  v_{n,\d})-a(x/|x|) |x|^{-2s} \eta_R  v_{n,\d} +  \eta_R  v_{n,\d}\Big) \,dx=0 .
\end{equation}
We claim that, for $n$ and $\delta$ fixed,
\begin{equation}\label{eq:claim_L2boundR}
\Ds( \eta_R v_{n,\d})\to \Ds v_{n,\d}\quad \text{in
}L^2(\R^N)\text{ as }R\to+\infty.
\end{equation}
By direct computations, we have
\begin{align*}
\Ds(\eta_R  v_{n,\d})(x)-\Ds v_{n,\d}(x)&= v_{n,\d}(x)\Ds\eta_R(x)  +(\eta_R(x)-1)  \Ds
v_{n,\d}(x)\\
&\qquad-c_{N,s}PV\int_{\R^N}\frac{(v_{n,\d}(x)-v_{n,\d}(y))(\eta_R(x)-\eta_R(y))}{|x-y|^{N+2s}}dy.
\end{align*}
Therefore, by H\"{o}lder's inequality, we get
\begin{align*}
\|\Ds(\eta_R  v_{n,\d})& -\Ds v_{n,\d}\|_{L^2(\R^N)}\\
&\leq \|v_{n,\d}\Ds\eta_R
\|_{L^2(\R^N)}+ \| (\eta_R-1)  \Ds v_{n,\d} \|_{L^2(\R^N)}\\
& + c_{N,s} \left( \int_{\R^N}\bigg(\,
\int_{\R^N}\frac{(v_{n,\d}(x)-v_{n,\d}(y))^2}{|x-y|^{N+2s}}dy\bigg)
\bigg(\int_{\R^N}\frac{(\eta_R(x)-\eta_R(y))^2}{|x-y|^{N+2s}}dy\bigg)
\,dx\right)^{1/2}.
\end{align*}
By scaling, we have, for some positive $C>0$ independent of $R$,
\begin{equation}\label{eq:DseR}
|\Ds\eta_R(x)|\leq C R^{-2s}\quad \text{for all }x\in \R^N.
\end{equation}
Next, we note that
\begin{align}\label{eq:Ds2}
c_{N,s}\int_{\R^N}\frac{(\eta_R(x)-\eta_R(y))^2}{|x-y|^{N+2s}}dy=-\Ds
(\eta_R^2)(x)+2\eta_R(x)\Ds\eta_R(x).
\end{align}
This implies, as above, that
$$
\int_{\R^N}\frac{(\eta_R(x)-\eta_R(y))^2}{|x-y|^{N+2s}}dy\leq C
R^{-2s},
$$
for some $C>0$. Therefore, using the above estimate and \eqref{eq:DseR}, we get
\begin{multline*}
\|\Ds(\eta_R  v_{n,\d})  -\Ds v_{n,\d}\|_{L^2(\R^N)}\\\leq  C (R^{-2s} +R^{-s}) \|v_{n,\delta}
\|_{H^s(\R^N)} + \|  (\eta_R-1) \Ds v_{n,\d} \|_{L^2(\R^N)}.
\end{multline*}
Hence \eqref{eq:claim_L2boundR} is proved. It follows that we can take the limit as $R\to \infty $ in
\eqref{eq:vndR} and use the dominated convergence theorem to obtain
\begin{equation}\label{eq:vnd}
  \int_{\R^N}w\Big(\Ds((1-\eta_\d)  v_n)-a(x/|x|) |x|^{-2s}(1-\eta_\d)   v_n + (1-\eta_\d)   v_n\Big)
  \,dx=0.
\end{equation}

\medskip\noindent
{\bf Step 3.}  We claim that
\begin{equation}\label{eq:claim_L2boundd}
\{\Ds( (1-\eta_\d) v_{n})\}_{\d \in(0,1)}\text{ is bounded in
}L^2(\R^N).
\end{equation}
As above, we have
  \begin{align*} \|\Ds(1-\eta_\d)&
v_n)\|_{L^2(\R^N)}\leq \|v_n\Ds\eta_\d  \|_{L^2(\R^N)}+ \|    \Ds
v_n \|_{L^2(\R^N)}  \\
& \quad+c_{N,s} \left( \int_{\R^N}\bigg(
\int_{\R^N}\frac{(v_{n}(x)-v_{n}(y))^2}{|x-y|^{N+2s}}dy\bigg)\bigg(
\int_{\R^N}\frac{(\eta_\d(x)-\eta_\d(y))^2}{|x-y|^{N+2s}}dy
\bigg)dx\right)^{1/2}.
\end{align*}
Let us estimate the first term in the right hand side of
the above inequality. To estimate it uniformly in $\d$, we use
Lemma \ref{lem:FW}  to get
$$
| \Ds\eta_\d |\leq C_{\eta,N,s}\frac{\d^{N}}{\d^{N+2s}+|x|^{N+2s}}.
$$
 Then we have
\begin{align*}
 \|v_n\Ds\eta_\d  \|_{L^2(\R^N)}^2&\leq  C\int_{B_1}v_n^2(x)\frac{\d^{2N}}{(\d^{N+2s}+|x|^{N+2s})^2}dx + C \|v_n  \|_{L^2(\R^N\setminus B_1)}^2 \\
  &\leq C \int_{B_\d} {v_n^2}\frac{\d^{2N}}{(\d^{N+2s}+|x|^{N+2s})^2}dx+ C \int_{B_1\setminus B_\d} {v_n^2}\frac{\d^{2N}}{(\d^{N+2s}+|x|^{N+2s})^2}dx\\
  &\quad+C  \|v_n  \|_{L^2(\R^N )}^2 \\
  &\leq C \int_{B_\d} v_n^2 \d^{-4s}dx+C  \int_{B_1\setminus B_\d} v_n^2 |x|^{2N}|x|^{-2N-4s}dx+C  \|v_n  \|_{L^2(\R^N )}^2 \\
  &\leq C \int_{B_\d} |x|^{-4s} v_n^2dx+C  \int_{B_1\setminus B_\d}  |x|^{-4s} v_n^2 dx+C  \|v_n  \|_{L^2(\R^N )}^2\\
  &\leq C\int_{B_1} |x|^{-4s} v_n^2dx+ C  \|v_n  \|_{L^2(\R^N )}^2,
\end{align*}
where $C$  is a positive constant independent on $\d$ (varying from
line to line). Hence by
\eqref{eq:vxest} and \eqref{eq:muanHarn},  we obtain
\begin{equation}\label{eq:vnDsed}
 \|v_n\Ds\eta_\d  \|_{L^2(\R^N)}^2\leq C+ C\|v_n  \|_{L^2(\R^N )}^2.
\end{equation}
 In addition,
by integration by parts, we have\footnote{
It is worth justifying the passage from \eqref{eq:I1I2I3} to
\eqref{eq:foot}. If $u\in H^{2s}(\R^N)$ and $g\in {\mathcal S}(\R^N)$,
the space of Schwarz functions, then
\begin{equation}\label{eq:1foot}
c_{N,s}\int_{\R^N}\bigg(
\int_{\R^N}\frac{(u(x)-u(y))^2}{|x-y|^{N+2s}}dy\bigg)g(x)dx=-\int_{\R^N}u^2(x) (-\D)^{s}g(x)\,dx+2\int_{
\R^N}u(x)g(x)(\Ds u)(x)\,dx.
\end{equation}
Indeed, we can approximate $u$ by smooth functions $u_n=\rho_n*u$  by
convolution with the standard mollifiers and consider $u_{n,R}=u_n
\eta_R$ so that $u_{n,R}\to u$ in $H^{2s}(\R^N)$ as
$n,R\to+\infty$, see \eqref{eq:claim_L2boundd}. Since $\Ds(u_{n,R}^2)\in L^2(\R^N)$ by  Lemma \ref{lem:FW}
we can pass to the limit in
 \begin{align*}
 &c_{N,s}\int_{\R^N}\bigg( \int_{\R^N}\frac{(u_{n,R} (x)-u_{n,R}
   (y))^2}{|x-y|^{N+2s}}dy\bigg)g(x)dx\\
&=\int_{\R^N}\big[-(-\D)^{s}(u_{n,R}^2)(x)+2 u_{n,R} (x)(\Ds u_{n,R})(x)\big]
 g(x)\,dx\\
&=-\int_{\R^N}u_{n,R}^2(x)
 (-\D)^{s}g(x)\,dx+2\int_{ \R^N}u_{n,R} (x)g(x)(\Ds u_{n,R})(x)\,dx
\end{align*}
to obtain \eqref{eq:1foot}.  }
 \begin{align}
\label{eq:I1I2I3}
&c_{N,s}\int_{\R^N}\bigg(
\int_{\R^N}\frac{(v_{n}(x)-v_{n}(y))^2}{|x-y|^{N+2s}}dy\bigg)\bigg(
\int_{\R^N}\frac{(\eta_\d(x)-\eta_\d(y))^2}{|x-y|^{N+2s}}dy
\bigg)dx  \\
&\nonumber =\int_{\R^N}[-\Ds (v_n^2)(x)+2v_n(x)\Ds v_n(x)]
[-\Ds
(\eta_\d^2)(x)+2\eta_\d(x)\Ds\eta_\d(x)]  dx  \\
&\label{eq:foot}=\int_{\R^N}v_n^2 (-\D)^{2s}\eta_\d^2-2\int_{
\R^N}v_n^2\Ds(\eta_\d\Ds\eta_\d) \\
&\nonumber\quad+ \int_{\R^N} 2v_n(x)\Ds v_n(x)  [-\Ds
(\eta_\d^2)(x)+2\eta_\d(x)\Ds\eta_\d(x)]  dx \\
&\notag\qquad\leq I_1+I_2+I_3,
\end{align}
where
$$
I_1= \int_{\R^N}v_n^2| (-\D)^{2s}\eta_\d^2|dx,
$$
$$
I_2=-2\int_{ \R^N}v_n^2\Ds(\eta_\d\Ds\eta_\d)dx
$$
and
$$
I_3= \left| \int_{\R^N} 2v_n(x)\Ds v_n(x)  [-\Ds
(\eta_\d^2)(x)+2\eta_\d(x)\Ds\eta_\d(x)]  dx\right|.
$$
 Now, for the last integral, we can use similar techniques as above
 to get
 \begin{align*}
I_3&\leq C \int_{\R^N}(|x|^{-2s}v_n+v_n+f)v_n
\frac{\d^{N}}{\d^{N+2s}+|x|^{N+2s}}dx\\
&\leq C \int_{B_1}|x|^{-4s}v_n^2dx+ \int_{\R^N\setminus B_1}v_n^2dx
+ C\int_{\R^N}fv_ndx,
\end{align*}
where $C>0$ (varying from line to line) is independent of $\delta$
(but could depend on $n$ and $\mathop{\rm Supp}f$).
Hence by \eqref{eq:vxest} and \eqref{eq:muanHarn}, we deduce that
\begin{equation}\label{eq:I3}
I_3\leq  C+ \int_{\R^N\setminus B_1}v_n^2dx + C \int_{\R^N}fv_ndx .
\end{equation}
We observe that $\Ds\eta_\d=\d^{-2s}\Ds \eta(\cdot/\d)\in
C^\infty(\R^N)$ so that $\eta_\d\Ds\eta_\d\in C^\infty_c(\R^N)$ and $\Ds(\eta_\d\Ds\eta_\d )=\d^{-4s} \Ds(\eta\Ds\eta
)(\cdot/\d). $ This implies that
 \begin{align*}
I_2&=-2\int_{ \R^N}v_n^2\Ds(\eta_\d\Ds\eta_\d)dx\leq C \int_{
\R^N}v_n^2\frac{\d^{N-2s}}{\d^{N+2s}+|x|^{N+2s}}dx\\
&\leq  C \int_{B_1}|x|^{-4s}v_n^2dx+ \int_{\R^N\setminus B_1}v_n^2dx
\end{align*}
and thus by \eqref{eq:vxest} and \eqref{eq:muanHarn}, we get
\begin{equation}\label{eq:I2}
I_2\leq  C+ C\int_{\R^N\setminus B_1}v_n^2dx   .
\end{equation}
Next, we estimate $I_1$. If $2s=1$ then
$$
I_1= \int_{\R^N}v_n^2| (-\D)^{2s}\eta_\d^2|dx\leq C \d^{-2}\int_{\d
\leq |x|\leq 2\d}v_n^2dx\leq C \int_{\R^N}|x|^{-2}v_n^2dx.
$$
If $2s<1$ then (using again Lemma \ref{lem:FW})
$$
I_1\leq C \int_{ \R^N}v_n^2\frac{\d^{N}}{\d^{N+4s}+|x|^{N+4s}}dx\leq
C \int_{B_1}|x|^{-4s}v_n^2dx+C \int_{\R^N\setminus B_1}v_n^2dx.
$$
Hence by \eqref{eq:vxest} and \eqref{eq:muanHarn}
\begin{equation}
I_1\leq  C+ C\int_{\R^N\setminus B_1}v_n^2dx   .
\end{equation}
If $2s>1$ then $0<2s-1<1$ so  $ (-\D)^{2s}\eta_\d^2=-(-\D)^{2s-1}(\D
\eta_\d^2 )$ which implies (see Lemma \ref{lem:FW}) that
$$
|(-\D)^{2s-1}(-\D \eta_\d^2 )|\leq C \d^{-2}
\frac{\d^{N}}{\d^{N+2(2s-1)}+|x|^{N+2(2s-1)}}.
$$
We then have, by similar estimates as a above,
$$
I_1\leq C \int_{ \R^N}v_n^2
\frac{\d^{N-2}}{\d^{N+4s-2}+|x|^{N+4s-2}}dx\leq C
\int_{B_1}|x|^{-4s}v_n^2dx+C \int_{\R^N\setminus B_1}v_n^2dx
$$
so that,  by \eqref{eq:vxest} and \eqref{eq:muanHarn}, we have
\begin{equation}
I_1\leq  C+ C\int_{\R^N\setminus B_1}v_n^2dx
\end{equation}
for $2s>1 $. We thus  conclude that for all $s\in (0,1)$
\begin{equation}\label{eq:I1}
I_1\leq  C+C \int_{\R^N\setminus B_1}v_n^2dx.
\end{equation}
Using the estimates  \eqref{eq:I1}, \eqref{eq:I2} and \eqref{eq:I3}
in \eqref{eq:I1I2I3}, together with \eqref{eq:vnDsed}, we get
 \eqref{eq:claim_L2boundd} as claimed.\\

\medskip\noindent
{\bf Step 4.} From \eqref{eq:claim_L2boundd} it follows that
$\Ds((1-\eta_\d)  v_n) \rightharpoonup \Ds v_n$ weakly in $L^2(\R^N)$
as $\delta \to 0^+$ (for any $n$ fixed).
Passing to the limit as $\delta\to0^+$ in \eqref{eq:vnd}, we then obtain, from the Dominated Convergence Theorem,
\eqref{eq:vxest} and \eqref{eq:muanHarn}, that
$$
  \int_{\R^N}w\Big(\Ds  v_n-a(x/|x|) |x|^{-2s}  v_n +   v_n\Big) \,dx=0 .
$$
Therefore, recalling \eqref{eq:vstf},
\be\label{eq:fweq}
\int_{\R^N}wf\, dx=\int_{\R^N}(a(x/|x|)-a_n(x/|x|) )|x|^{-2s}  v_nw\,dx.
\ee
By H\"{o}lder's inequality, Fubini's theorem, and estimates
\eqref{eq:vxest} and \eqref{eq:vxestinfty}
\begin{align*}
&\left|  \int_{\R^N}\frac{a(x/|x|)-a_n(x/|x|)}{ |x|^{2s}}  v_nw
\right|^2\leq \|w\|_{L^2(\R^N)}^2 \Big\|\frac{a(x/|x|)-a_n(x/|x|)}{ |x|^{2s}}
v_n\Big\|_{L^2(\R^N)}^2\\
&\quad\leq \|w\|_{L^2(\R^N)}^2\int_0^\infty r^{-4s+N-1}\left(
\int_{\S^{N-1}}|v_n(r\th')|^{{\frac{2q}{q-2}}}dS' \right)^{\!\!\frac{q-2}q}\left(
\int_{\S^{N-1}}|a_n-a|^{q} dS' \right)^{\!\!\frac2q}dr\\
&\quad \leq \|w\|_{L^2(\R^N)}^2C_0^2|\S^{N-1}|^{(q-2)/q}\| a_n-a\|_{L^{q }(\S^{N-1})}^2 \int_0^{r_0}
r^{-4s+N-1+2\g_n}dr \\
& \quad+
 \|w\|_{L^2(\R^N)}^2C_1^2|\S^{N-1}|^{(q-2)/q}\| a_n-a\|_{L^{q }(\S^{N-1})}^2 \int_{r_0}^\infty
r^{-4s+N-1+2\a_n}dr\\
&\quad \leq   C\| a_n-a\|_{L^{q }(\S^{N-1})}^2.
\end{align*}
From this,  \eqref{eq:fweq} and \eqref{eq:antoa}, we deduce that
$$
\int_{\R^N}fw\, dx=0.
$$
We have then proved that $\int_{\R^N}fw\, dx =0$ for every $f\in C^\infty_c(\R^N\setminus\{0\})$, $f\geq0$,
$f\neq0$. This implies that $w\equiv 0$ which  leads  to a contradiction.

\smallskip
\noindent
 \textbf{Case 2: } \be\label{eq:musHmeqNot}
 -\mu_1(a)=
\frac{(N-2s)^2}{4}-s^2.
  \ee
As in Case 1, we argue by contradiction and assume that $A$ is not
essentially self-adjoint; as observed above, this implies that there
exists $w\in L^2(\R^N)\setminus\{0\}$ such that
$$\Ds w-a(x/|x|) |x|^{-2s}w +w =0$$  in the sense of distributions in $\R^N\setminus\{0\}$.

Let  $a_n\in C^\infty(\S^{N-1})$ be as in \eqref{eq:antoa}. Let $\s\in(0,1)$ and notice that $ \mu_1(a-\s)>\mu_1(a)$  so
that
 $$-\mu_1(a-\s)<\frac{(N-2s)^2}{4}-s^2.$$
 Given $f\in C^\infty_c(\R^N\setminus\{0 \})$ such that $f\geq0$, $f\not\equiv0$, by the Lax-Milgram
theorem, for every $n\in\N$ and $\sigma\in(0,1)$, there exists
$v^\s_n\in H^s(\R^N)$ ($v^\s_n>0$ in $\R^N\setminus\{0\}$) solution to
$$
  \Ds v^\s_n -  |x|^{-2s}(a_n(x/|x|)-\s)v^\s_n+v^\s_n =f.
$$
It is then not difficult to check that $(v_n^\s)_n$ is
bounded in $H^s(\R^N)$ and converges (weakly in $H^s(\R^N)$) to some
$v^\s$ weakly solving
$$
  \Ds v^\s -  |x|^{-2s}(a(x/|x|)-\s)v^\s+v^\s =f.
$$
Arguing as in Lemma \ref{l:stima}, we have that there exist $r_1,C_2>0$  independent on $\s$ and
$n$ such that
\be\label{eq:Vdest}
v^\s_n(x)|\leq C_2 |x|^{\g_n(\s)},\quad v^\s(x)\leq C_2
|x|^{\g(\s)} \quad\text{for all }x\in B_{r_1}\setminus\{0\},
 \ee
 where
\[
\g_n(\s)= -\tfrac{N-2s}{2}+\sqrt{\big(\tfrac{N-2s}
    {2}\big)^{2}+\mu_{1}(a_n-\sigma)},\quad
\g(\s)= -\tfrac{N-2s}{2}+\sqrt{\big(\tfrac{N-2s}
    {2}\big)^{2}+\mu_{1}(a-\sigma)}.
\]
Notice that, provided $n$ is large,
$$-s+\sqrt{\bigg(\frac{N-2s}
    {2}\bigg)^{2}+\mu_{1}(a_n-\s)}>\e_0>0, $$
for some $\e_0>0$ depending on $\sigma$ (but independent of $n$).
 Therefore, using
similar arguments as in Case 1, we get
$$
\left| \int_{\R^N}f w dx \right|\leq \s \int_{\R^N} |x|^{-2s}|w|
v^\s dx .
$$
Hence, by  H\"{o}lder's inequality
$$
\left| \int_{\R^N}f w \, dx \right|\leq \|w\|_{L^2(\R^N)} \s  \||x|^{-2s}
v^\s\|_{L^2(\R^N)}.
$$
Since $v^\s$ is bounded in $H^s(\R^N)$, so $  \||x|^{-2s}  v^\s
\|_{L^2(\R^N\setminus B_1)}$ can be uniformly bounded in $\s$, we
infer that
$$
\left| \int_{\R^N}f w \,dx \right|\leq
\|w\|_{L^2(\R^N)} \s C_{a,N,s,f}+
\|w\|_{L^2(\R^N)} \s
\||x|^{-2s}
v^\s\|_{L^2(B_1)}
$$
for all $\sigma\in(0,1)$.
Using \eqref{eq:Vdest} we deduce that, for all $\sigma\in(0,1)$ and for some $C>0$ independent of $\sigma$,
\begin{equation}\label{eq:8}
\left| \int_{\R^N}f w \,dx \right|\leq C \frac{\s}{\left(
-s+\sqrt{\frac{(N-2s)^2}{4}+\mu_1(a -\s)} \right)^{1/2}} + C \s .
\end{equation}
Using the variational characterization of  $\mu_1(a-\s) $ and
$\mu_1(a)$ (see also \cite[Proof of Lemma 2.1]{FF1} for convergences
of related eigenfunctions), we deduce that
$$
 c_2\s+\mu_1(a)\leq  \mu_1(a-\s)\leq \mu_1(a)+c_1\s,
$$
where  $c_2,c_1$ are positive constants independent on $\s$. Passing
to the limit in \eqref{eq:8} as $\s\to 0^+$ yields that $\int_{\R^N}fwdx=0$. This
then implies that $w=0$ thus giving rise to a contradiction. \QED

\medskip\noindent  {\bf Proof of Theorem \ref{t:esa}: sufficiency of condition \eqref{eq:11}.}
  In view of Theorem \ref{th:essself}, condition \eqref{eq:11} implies that $A=\Ds-a(x/|x|)
|x|^{-2s}$ with domain $C^\infty_c(\R^N\setminus\{0\})$ is
essentially self-adjoint in $L^2(\R^N)$. Then for every $m\geq 0$ also
the operator $\Dsm-a(x/|x|)
|x|^{-2s}$ is essentially self-adjoint in $L^2(\R^N)$ by Lemma
\ref{eq:meq0}, part $(i)$. \QED

\begin{Remark}\label{rem:self-h}
If  $h\in L^\infty_{loc}(\R^N)\cap L^p(B_r)$, for some $p >N/(2s)$, $r>0$  and $h$ is bounded in a
neighborhood of $\infty$, then the arguments proving Theorem
\ref{th:essself}  above can be adapted to
prove that if $-\mu_1(a)\leq
\frac{(N-2s)^2}{4}-s^2$ and  $N>2s$ then the operator $\Dsm-a(x/|x|)
 |x|^{-2s}+h(x)$ with domain $C^\infty_c(\R^N\setminus\{0\})$ is
essentially self-adjoint in $L^2(\R^N)$. Indeed the estimate of Lemma
\ref{l:stima} still holds for the operator
$\Ds-a(x/|x|)|x|^{-2s}+h$, see  \cite[Lemma 5.11]{FF1}.
\end{Remark}

\section{ Non-essential self-adjointness} \label{s:non-ess}

The following lemma will be crucial in our proof.
\begin{Lemma}\label{lem:Dsbf}
Let $a\in L^\infty(\S^{N-1})$ and  $b>0$. Let $\psi_1$ be an
eigenfunction of problem \eqref{eq:1eig} corresponding to first
eigenvalue $\mu_1(a)$ in \eqref{firsteig}.
Let $\nu_1={\sqrt{(N-2s)^2/4+\mu_1(a) }}$ and assume that
\[
-\mu_1(a)>
\frac{(N-2s)^2}{4}-s^2,\quad\text{i.e. }s>\nu_1.
\]
 For $z=(t,x)\in \RNp$ define $f(z) = \psi_1(z/|z|) |z|^{\frac{2s-N}{2}}K_{\nu_1
 }(\sqrt{b}|z|)$.  Then
\be\label{eq:fundHard}
\begin{cases}
-\mathrm{div}(t^{1-2s}\nabla f)+ t^{1-2s} b f= 0,   &\mbox{in } \RNp, \\
-\dlim_{t\rightarrow 0^+}t^{1-2s}\pa_tf= \k_s
a(x/|x|)|x|^{-2s}f, &\mbox{on } \R^N\setminus\{0\}.
\end{cases}
\ee
 In addition
\be\label{eq:Dsbf}
(-\D+b)^sf- a(x/|x|)|x|^{-2s}f=0 \quad \mbox{ in  }
\calD'(\R^N\setminus\{0\}),
\ee
i.e.
\[
\int_{\R^N}f(0,x)\bigg((-\D+b)^s\varphi(x)-
\frac{a(x/|x|)}{|x|^{2s}}\varphi(x)\bigg)\,dx=0\quad\text{for all }\varphi \in
C^{\infty}_c(\R^N\setminus\{0\}).
\]
\end{Lemma}

\begin{Remark}
We notice that  the conclusion of the above lemma might not be true
if $b=0$. Here, we have the property that $K_{\nu_1}$ decays
exponentially at infinity which plays a central role in the proof.
\end{Remark}
\proof Direct computations  using polar coordinates, see \cite{FF},
prove the first assertion.

Let $\vp\in C^\infty_c(\R^N\setminus\{0\})$.  We now consider the extension
$\Phi(t,x)=(P_{\sqrt{b}}(t,\cdot)*\vp )(x)$, where $P_{\sqrt{b}}$
is the Bessel Kernel, see Section \ref{s:sp}. We have that
$\Phi\in H^1(\RNp; t^{1-2s})$ and
 \be\label{eq:Phin}
\begin{cases}
-\mathrm{div}(t^{1-2s}\nabla \Phi)+ t^{1-2s} b \Phi= 0, &\mbox{in } \RNp, \\
\Phi=\phi,  &\mbox{on } \R^N, \\
-\dlim_{t\rightarrow 0^+}t^{1-2s}\pa_t\Phi = \k_s (-\D+b)^s
\phi, &\mbox{on } \R^N.
\end{cases}
\ee
We multiply the first equation of \eqref{eq:Phin} by $f$ and
integrate  by parts over $\R^N\times(\rho,\infty)$ for $\rho>0$ to
get
\begin{align}\label{eq:estnVPnef}
2\int_{\R^N} \rho^{1-2s}
 \pa_t\Phi(\rho,x) f(\rho,x) dx
 - 2\int_{\R^N}  \rho^{1-2s} \pa_t f(\rho,x) \Phi(\rho,x) dx=0.
\end{align}
 By  \eqref{eq:decKnu0} and \eqref{eq:decKnuInf},
$|f(\rho,x)|\leq C |x|^{ \frac{2s-N}{2}-\nu_1}$  for all $x\in
B_1$ and  $|f(\rho,x)|\leq C e^{-\frac{\sqrt{b}}2|x|}$ for every
$x\in\R^N\setminus B_1$ and $\rho\in (0,1)$. Since $\frac{N+2s}{2}-\nu_1>0$,
$f(0,\cdot)\in L^1(B_1)$.  Since $\vp\in C^\infty_c(\R^N)$, we
have that $ (-\D+b)^s\vp\in C(\R^N)\cap L^\infty(\R^N)$. By using  \cite{FF1},
$$
|\rho^{1-2s}\de_t\Phi(\rho,x)|\leq C,
 \quad \text{for all }\rho\in(0,1),\
x\in \R^N.
$$
Therefore  we can apply the dominated convergence theorem and use \eqref{eq:Phin} to get
\be\label{eq:limdePhinf}
\lim_{\rho\to 0}\int_{\R^N} \rho^{1-2s}
 \pa_t\Phi(\rho,x) f(\rho,x) dx= -\k_s\int_{\R^N}(-\D+b)^s\phi(x)f(0,x) dx.
\ee
It now remains to prove that
$$
\lim_{\rho\to 0} \int_{\R^N}  \rho^{1-2s} \pa_t f(\rho,x) \Phi(\rho,x) \, dx=\k_s \int_{\R^N} a(x/|x|)|x|^{-2s}\phi(x)f(0,x)  \,dx
$$
which completes the proof. This will be done in the sequel.

By direct computations we have
\begin{align*}
 t^{1-2s}\de_t f(z)&=|z|^{\frac{2s-N}{2}-2s}K_{
\nu_1}(\sqrt{b}|z|)
\th_1^{1-2s}\n_{\S^N}\psi_1(z/|z|)\cdot  {\mathbf e}_1\\
&\quad +[\tfrac{2s-N}{2}
t^{2-2s}|z|^{\frac{2s-N}{2}-2}K_{\nu_1 }(\sqrt{b}|z|) +\sqrt{b}
t^{2-2s}|z|^{\frac{2s-N}{2}-1}K_{\nu_1 }'(\sqrt{b}|z|)]
\psi_1(z/|z|)\\
& := H(t,x)+J(t,x),
\end{align*}
where, for $z=(t,x)$, we define
$$
H(t,x):= |z|^{\frac{2s-N}{2}-2s}K_{ \nu_1}(\sqrt{b}|z|)
\th_1^{1-2s}\n_{\S^N}\psi_1(z/|z|)\cdot  {\mathbf e}_1
$$
and
$$
J(t,x):=\Big[\tfrac{2s-N}{2} t^{2-2s}|z|^{\frac{2s-N}{2}-2}K_{\nu_1
}(\sqrt{b}|z|) +\sqrt{b} t^{2-2s}|z|^{\frac{2s-N}{2}-1}K_{\nu_1
}'(\sqrt{b}|z|)\Big] \psi_1(z/|z|).
$$

 First we recall that $\psi_1\in  C^{0,\a}(\ov{\S^N_+})$, see
\cite{FF1}. Once again by \eqref{eq:decKnu0} and
\eqref{eq:decKnuInf} together with the fact that
$K_{\nu_1}'=-\frac{\nu_1}{r}K_{\nu_1}-K_{\nu_1-1}$, it is plain that
for $|x|\leq 1$
\begin{equation}\label{eq:J}
|J(\rho,x )|\leq c \rho^{2-2s} |(\rho,x)|^{\frac{2s-N}{2}-2-\nu_1}
\end{equation}
while  for $|x|\geq 1$
$$
|J(\rho,x )|\leq c \rho^{2-2s}   e^{-\sqrt{b}/2|x|}.
$$
 It is then clear that
$$
\lim_{\rho\to 0}\int_{\R^N\setminus B_1}    J(\rho,x) \Phi(\rho,x)
dx=0.
$$
Recalling the notations in Section \ref{s:sp},    for $\rho\leq 1$,
we claim that
\be\label{eq:taylorsBess}
\Phi(\rho,x)=\phi(x)\vartheta(\sqrt{b}\rho)+O(\rho^{2s}),\quad
\text{for all } x\in \R^N.
\ee
To see this, we use change of variables and  \eqref{eq:intPm}, to
get, up some normalization constant,
\begin{align*}
&\Phi(\rho,x)=( P_{\sqrt{b}}(\rho,\cdot)*\vp) (x)=
\int_{\R^N}P_{\sqrt{b}}(\rho,y)\vp(x+y)dy=\vp(x)\int_{\R^N}P_{\sqrt{b}}(\rho,y)dy\\
&+ \int_{\R^N}  P_{\sqrt{b}}(\rho,y)\n \vp(x)\cdot
y dy + \int_{\R^N}\int_{0}^{1}\int_{0}^{1}\t
P_{\sqrt{b}}(\rho,y)D^2\vp(x+\t r y)[y,y]  d\t dr
dy\\
\qquad &=\vp(x)\vartheta(\sqrt{b}\rho)+0+
\int_{\R^N}\int_{0}^{1}\int_{0}^{1}\t  P_{\sqrt{b}}(\rho,y)D^2\vp(x+\t r
y)[y,y] d\t dr dy.
\end{align*}
Therefore
$$
|\Phi(\rho,x)- \vp(x)\vartheta(\sqrt{b}\rho)|\leq C \rho^{2s}
  \int_{|y|\leq 1}  |y|^{-N-2s+2} dy+ C \rho^{2s}
  \int_{|y|\geq 1}  e^{-\sqrt{b}/2|y|}  dy,
$$
 thus \eqref{eq:taylorsBess} is proved. From  \eqref{eq:taylorsBess}  together with \eqref{eq:J}, we
 deduce that
$$
 |J(\rho,x )\Phi(\rho,x)|\leq C |
\phi(x)| |x|^{\frac{2s-N}{2}-2-\nu_1}     + C
|(\rho,x)|^{\frac{2s-N}{2}-\nu_1} \quad \text{for all } x \in B_1.
$$
Therefore
$$
 |J(\rho,x )\Phi(\rho,x)|\leq
 C +C|x|^{\frac{2s-N}{2}-\nu_1}\quad \text{for all } x \in B_1 .
$$
 The dominated convergence theorem then implies that
$$
\lim_{\rho\to 0} \int_{B_1} J(\rho,x )\Phi(\rho,x)dx =0.
$$
Hence
\[
\lim_{\rho\to 0}\int_{\R^N }    J(\rho,x) \Phi(\rho,x) dx=0.
\]
It remains now to pass the limit as $\rho\to 0$ in the integral
$\int_{\R^N}H(\rho,x)\Phi(\rho,x)dx$. To this end, we first claim that
 \be \label{eq:claim}
\th_1^{1-2s}\n_{\S^N}\psi_1\cdot {\mathbf e}_1
\in
L^\infty( \ov{\S^N_+}).
\ee
 To prove this claim, we consider
\[
 g(z) = \psi_1(z/|z|)  |z|^{\frac{2s-N}{2}}I_{{\sqrt{
(N-2s)^2/4+\mu_1(a) }} }(\sqrt{b}|z|)
\]
which  satisfies
\be\label{eq:gundHard}
\begin{cases}
-\mathrm{div}(t^{1-2s}\nabla g)+ t^{1-2s} b\,  g= 0, &\mbox{in } \RNp, \\
-\dlim_{t\rightarrow 0^+}t^{1-2s}\pa_tg= \k_s
a(x/|x|)|x|^{-2s}g, &\mbox{on } \R^N\setminus\{0\},
\end{cases}
\ee where $I$ is the modified Bessel function of first kind. Using
its decay property near the origin (see \cite{El}), we see that
$$
g(z)\leq C |z|^{\frac{2s-N}{2}+\nu_1}, \quad \text{for all } z\in
\ov{B_2^+},
$$
 implying  $ |z|^{-1}g\in
L^2(B_2^+; t^{1-2s})$ and $ |x|^{-s}g(0,\cdot)\in L^2(B_2)$.
Using standard integration by parts, we can deduce that
$g\in H^1(B_2^+;t^{1-2s})$. Then by \cite{FF1} it follows
that $t^{1-2s}\de_tg\in  L^\infty(\ov{B_{3/2}^+\setminus
B^+_{1/2}})$. As above, by direct computations,
\begin{align*}
t^{1-2s}\de_t g(z)&=|z|^{\frac{2s-N}{2}-2s}I_{ \nu_1}(\sqrt{b}|z|)
\th_1^{1-2s}\n_{\S^N}\psi_1(z/|z|)\cdot  {\mathbf e}_1\\
& \qquad+\big[\tfrac{2s-N}{2} t^{2-2s}|z|^{\frac{2s-N}{2}-2}I_{\nu_1
}(\sqrt{b}|z|) + \sqrt{b}\,t^{2-2s}|z|^{\frac{2s-N}{2}-1}I_{\nu_1
}'(\sqrt{b}|z|)\big] \psi_1(z/|z|) .
\end{align*}
Evaluating at $|z|= 1$ and using the fact that $I_{\nu_1
}(\sqrt{b})\neq 0$, we see that
$$|\th_1^{1-2s}\n_{\S^N}\psi_1(z)\cdot
 {\mathbf e}_1|\leq C( \| \psi_1\|_ {L^\infty(\S^N_+ )}+1), \quad \text{for all } z\in
\ov{\S^N_+},
$$
and  claim  \eqref{eq:claim} is proved.\\

\noindent
 Using the Taylor expansion \eqref{eq:taylorsBess}, and similar arguments as
 above, we obtain
 \begin{align*}
\lim_{\rho\to 0}\int_{\R^N}H(\rho,x)\Phi(\rho,x)dx&=
-\kappa_s\int_{\R^N}|x|^{\frac{2s-N}{2}-2s}K_{ \nu_1}(\sqrt{b}|x|)
a(x/|x|)\psi_1(x/|x|)\phi(x)dx\\
&= -\kappa_s\int_{\R^N}a(x/|x|)|x|^{ -2s}f(0,x)\phi(x)dx.
 \end{align*}
Using this, \eqref{eq:limdePhinf} and \eqref{eq:estnVPnef}, we get
$$
\int_{\R^N}a(x/|x|)|x|^{
-2s}f(0,x)\phi(x)dx=\int_{\R^N}(-\D+b)^s\phi(x)f(0,x) dx
$$
which is \eqref{eq:Dsbf}. \QED

\begin{Theorem}[\bf Necessity of condition
  \eqref{eq:11} of Theorem \ref{t:esa}]\label{t:nonsa}
Let $N>2s$, $m\geq0$, and $a\in L^\infty(\S^{N-1})$.
  Then the operator
 $A'=\Dsm-a(x/|x|) |x|^{-2s}$ with domain $C^\infty_c(\R^N\setminus\{0\})$ is not
essentially self-adjoint in $L^2(\R^N)$ if
\be\label{eq:assumpnS}
 -\mu_1(a)>
\frac{(N-2s)^2}{4}-s^2.
\ee
\end{Theorem}
\proof Assume by contradiction that \eqref{eq:assumpnS} holds and $ A'$ is essentially
self-adjoint. By part $(i)$  of Lemma \ref{eq:meq0} also $(-\D +b)^s
-a(x/|x|) |x|^{-2s}$ is  essentially
self-adjoint in $L^2(\R^N)$ for every $b>0$; then by assumption
\eqref{eq:pot}, Remark \ref{sec:essent-self-adjo-1},  and part $(ii)$ of
Lemma \ref{eq:meq0}, we have that the operator $ (-\D +b)^s -a(x/|x|) |x|^{-2s}$ has dense range in $L^2(\R^N)$.

Let $\psi_1$ be a positive eigenfunction of problem \eqref{eq:1eig}
associated to the first eigenvalue $\mu_1(a)$ defined in
\eqref{firsteig}.  For $z=(t,x)$ let
\[
 f(z) = \psi_1(z/|z|) |z|^{\frac{2s-N}{2}}K_{{\sqrt{
      (N-2s)^2/{4}+\mu_1(a) }} }(\sqrt{b}\,|z|).
\]
We observe that by
\eqref{eq:assumpnS}, \eqref{eq:decKnu0} and \eqref{eq:decKnuInf}
\be
\label{eq:dhfL2} f(0,\cdot)\in L^2(\R^N) .
\ee
Since $ (-\D +b)^s -a(x/|x|) |x|^{-2s}$ has dense range in
$L^2(\R^N)$ (as assumed above for the contradiction), there exists
$\vp_n\in C^\infty_c(\R^N\setminus\{0\})$ such that
$(-\D +b)^s \vp_n-a(x/|x|) |x|^{-2s}\vp_n \to f(0,\cdot)$ in $L^2(\R^N)$, so that for every
$\e>0$ there exists $n(\e)$ such that
$$
\|(-\D +b)^s \vp_n-a(x/|x|) |x|^{-2s}\vp_n - f(0,\cdot)\|_{L^2(\R^N)}<\e\quad
\text{for every } n\geq n(\e).
$$
This implies that
\begin{equation*}
-2\int_{\R^N}\Big((-\D +b)^s \vp_n (x) -a(x/|x|)|x|^{-2s} \vp_n(x)
\Big)f(0,x)\,dx+\|
f\|_{L^2(\R^N)}^2 <\e^2
\end{equation*}
{ for every } $n\geq n(\e)$.
By \eqref{eq:Dsbf} we obtain that for every $\e>0$
\begin{equation*}
\| f\|_{L^2(\R^N)}^2  <\e^2
\end{equation*}
which is impossible.
 \QED

    \label{References}

\end{document}